\documentclass[a4,12pt]{article}
\usepackage{geometry,times}
\usepackage{graphicx}
\usepackage{amssymb}
\usepackage{epstopdf}
\usepackage{makeidx}
\usepackage{amsmath}
\usepackage{amsthm}
\usepackage[T1]{fontenc}
\usepackage{amsfonts}
\usepackage{hyperref}
\usepackage{marvosym}
\usepackage{bbm}

\theoremstyle{definition}

\newcommand{\IE}{\mathbb{E}}

\newcommand{\IP}{\mathbb{P}}

\newcommand{\ID}{\mathbb{D}}
\newcommand{\IV}{\mathbb{V}}

\renewcommand{\P}{\mathcal{P}}

\newcommand{\tn}{\textnormal}
\newcommand{\Xn}{X_{\nu\;\! i}}
\makeatletter

\renewenvironment{thebibliography}[1]
      {\section*{\refname
         \@mkboth{\refname}{\refname}}\small
       \list{\@biblabel{\@arabic\c@enumiv}}%
            {\settowidth\labelwidth{\@biblabel{#1}}%
             \leftmargin\labelwidth
             \itemsep0pt\parsep.5ex
             \advance\leftmargin\labelsep
             \@openbib@code
             \usecounter{enumiv}%
             \let\p@enumiv\@empty
             \renewcommand\theenumiv{\@arabic\c@enumiv}}%
       \sloppy\clubpenalty4000\widowpenalty4000%
       \sfcode`\.\@m}
      {\def\@noitemerr
        {\@latex@warning{Empty `thebibliography' environment}}%
       \endlist}

\makeatother

\title{The fate of the square root law for correlated voting}

\author{Werner Kirsch\footnote{Werner Kirsch, \Letter Fakult\"at f\"ur Mathematik und Informatik, FernUniversit\"at Hagen, D-58095 Hagen, Germany, werner.kirsch@fernuni-hagen.de}\,
$\cdot$ Jessica Langner\footnote{Jessica Langner, \Letter Fakult\"at f\"ur Mathematik und Informatik, FernUniversit\"at Hagen, D-58095 Hagen, Germany, jessica.langner@fernuni-hagen.de}}

\begin{document}
\maketitle

\begin{abstract}
We consider two-tier voting system and try to determine optimal weights for a fair representation in such systems. A prominent example of such a voting system is the Council of Ministers of the European Union. Under the assumption of independence of the voters, the square root law gives a fair distribution of power (based on the Penrose-Banzhaf power index) and a fair distribution of weights (based on the concept of the majority deficit), both given in the book by Felsenthal and Machover.\\
In this paper, special emphasis is given to the case of correlated voters. The cooperative behaviour of the voters is modeled by suitable adoptions of spin systems known from statistical physics. Under certain assumptions we are able to compute the optimal weights as well as the average deviation of the council's vote from the public vote which we call the democracy deficit.
\end{abstract}

\paragraph*{Acknowledgement} This paper has been presented at the
Leverhulme Trust sponsored Voting Power in Practice Symposium held at the London School
of Economics, 20--22 March 2011.


\section{Introduction}

In this paper, we consider two-tier voting systems. The first level of such a systems usually consists of the voters in a country or an association of countries.
The voters in each constituency (or member country) are represented by a delegate in the second level voting system, the council. Delegates
in the council are given a voting weight which as a rule depends on the population of the constituency they represent.

Examples of such two-tier voting systems are the Council of Ministers of the European Union, the Electoral College in the USA and the `Bundesrat', the state chamber of Germany's parliamentary system. In each case we assume that the representatives vote according to the majority vote in their respective constituency.

What is a fair voting weight for a delegate in a council? This question arises immediately in all these examples. It seems self-evident that for a fair voting system the voting outcome in the council should agree with the result of a popular vote. The US presidential elections 2000 show that this is not always the case. While Al Gore won the public vote the majority in the Electoral College elected Georg W. Bush as the 43rd president of the USA. The difference between the voting result in the council and the public vote is called the `democracy deficit'.

In fact, it is not hard to see, that \emph{no} voting system for the council can \emph{guarantee} that the vote in the council and the public vote agree. In other words, no matter how we choose the voting weights for the council members, the democracy deficit cannot be zero for \emph{all} possible distributions of `yes'- and `no'- votes among the voters. Thus, the best one can do is to minimize the expected democracy deficit, i.~e., the difference between the vote in the council and popular vote. Obviously, the term `expected' needs a careful interpretation. If one assumes that all voters cast their votes independently of each other then one can show that the expected democracy deficit is minimized if the voting weight of a representative is chosen proportional to the \emph{square root} $\sqrt{N_\nu}$ of the population ($N_\nu$) of the respective country (with number $\nu$).

This is (one version of) the celebrated `square root law' by Penrose (see \\\cite{FeMa1998} and \cite{Penrose1946}). In this paper, we go beyond the square root law by dropping the assumption of the voters' independence. We apply two different schemes to model the correlation between the voters. In our main model we assume that the voters are influenced by a `common belief' of the society or -which is the same, technically speaking- by a strong group of opinion makers. We call this system the CBM (for `common belief model' or `collective bias model') (see \cite{Kirsch2007}). The CBM can be looked upon as a generalization of a model proposed by Straffin \cite{Straffin1977} in connection with the Shapley-Shubik power index (see \cite{ShSh1954}). The other model we look at takes into account that voters influence each other. It is based on a model (the Curie-Weiss Model) for ferromagnetic behaviour taken from statistical physics (see \cite{Kirsch2007} and cf. \cite{Ellis1985,Thompson1972}).

If we assume that the voters in different countries vote independently of each other, we can compute the optimal voting weights in terms of the expected margins of the voting outcome in the countries. For the CBM the optimal weights are proportional to the population $N_\nu$. We also compute the expected democracy deficit for these models (for large $N_\nu$).

Under the assumption that the voters influence each other also across country borders (according to the CBM) we can also compute the expected democracy deficit asymptotically. It turns out that in this case any voting weight is as good as any other one. In other words, on an asymptotical scale \emph{any} distribution of voting weights is close to optimal.

\section{The General Model}

We consider a situation where $M$ states (countries, constituencies) form a federation. The states are labeled by Greek characters, e.~g., $\nu, \kappa, \ldots$.
The number of voters (population) of the state $\nu$ is denoted by $N_\nu$. Consequently, the total population of the union is given by $N=\sum_{\nu=1}^M N_\nu$.

We represent the vote of the voter $i$ in state $\nu$ by $X_{\nu\,i}$. This voter may vote either `yes', in which case we set $X_{\nu\,i}=1$ or `no' encoded as
$X_{\nu\,i}=-1$. Consequently, the result of a simple majority voting in the state $\nu$ is represented by the sum $S_\nu=\sum_{i=1}^{N_\nu} X_{\nu\,i}$. A voting in that state is affirmative if $S_\nu>0$. For the simplicity of notation and to avoid nonsignificant technicalities we assume that all $N_\nu$ are odd numbers, this excludes a draw described by $S_\nu=0$.

We denote the voting decision in the state $\nu$ by $\chi_\nu=\chi_\nu(S_\nu)$ which we set equal to $1$ if $S_\nu>0$ and equal to $-1$ if $S_\nu\leq 0$. Thus, the representative of state $\nu$ will vote `yes' if $\chi_\nu=1$ and `no' if $\chi_\nu=-1$. For later use we note that $\chi_\nu S_\nu = |S_\nu|$.

If we denote the voting weight for state $\nu$ in the council by $g_\nu$ then the voting result in the council is given by
\begin{equation}\label{CouncVote}
C~=~\sum_{\nu=1}^M\;g_\nu\,\chi_\nu\,.
\end{equation}
This voting result has to be compared with the popular vote given by
\begin{equation}\label{PopVote}
P~=~\sum_{\nu=1}^M\;S_\nu\,.
\end{equation}
We call the absolute value of the difference between $C$ and $P$ the \emph{democracy deficit} and denote it by $\Delta$
\begin{align}\label{dd}
\Delta~&=~|\,C-P\,|\\
&=\Big|\,\sum_{\nu=1}^M\;g_\nu\,\chi_\nu\;-\;\sum_{\nu=1}^M\;S_\nu\,\Big|\,.
\end{align}

The democracy deficit $\Delta$ depends explicitly on the voting weights $g_1,\ldots,g_M$. The voting weights should be chosen in such a way that the democracy
deficit is as small as possible.

The voting results $\Xn$ are the voter's reaction on a particular proposal $\omega$. Hence, the democracy deficit $\Delta$ depends on the given proposal $\omega$
as well. It is easy to choose the weights $g_\nu$ such that $\Delta$ vanishes for a \emph{given} proposal. But our goal is to optimize the weights in such a way
that $\Delta$ is small for \emph{most} proposals. Thus, we look at the \emph{expected} value of $\Delta^2$, denoted by
\begin{equation}\label{expDelta}
\ID~:=~\IE\Big(\Delta^2\Big)\,.
\end{equation}
We will call $\ID$ the \emph{expected democracy deficit} in the following (instead of the correct but clumsy `expected square of the democracy deficit').

By looking at expectation values we regard the proposals as random input to the voting system. Hence the probability that the next proposal to the system is a particular proposal $\omega$ is determined by a probability rule. We assume that there is no bias to certain proposals, in particular any proposal and its counterproposal have the same probability.

The voting system reacts in a deterministic (and rational) way to this random input. The voting results as well as the democracy deficit are therefore (otherwise deterministic) functions of the random input, the proposal. The voting outcome is a vector in the space $\Omega=\{-1,1\}^N$, where $N$ is the total number of voters and the probability distribution of the proposals equips $\Omega$ with probability distribution $\IP$ as well, namely the probability  of a given outcome $(X_1,\ldots, X_N)$ is the probability of all proposals $\omega$ that lead to that outcome. Since the voters react rationally they vote $-1$ on the opposite to a proposal they would favour and vice versa. Hence the probability distribution $\IP$ satisfies

\begin{equation}
\IP(X_1,\ldots, X_N)~=~\IP(-X_1,\ldots, -X_N)\,.
\end{equation}
We call such a measure a \emph{voting measure}. For any voting measure we have $\IP(X_i=1)=\IP(X_i=-1)=\tfrac{1}{2}$, but probabilities concerning more than one voter, like $\IP(X_1=1\  \tn{and}\ X_2=1)$ cannot be computed from the mere assumption that $\IP$ is a voting measure. Such events concern the correlation structure of the measure and they have yet to be fixed depending on the situation at hand. One possible specification is the assumption that all voters act independently of each other.
This leads to the property that
\begin{align*}
\IP(X_1=1\  \tn{and}\ X_2=1)\,=\,\IP(X_1=1)\cdot\IP(X_2=1)\,=\,\tfrac{1}{4}\,.
\end{align*}
More generally, under the assumption of independence we have
\begin{equation}\label{indep}
\IP\big(X_1=\xi_1, X_2=\xi_2,\ldots, X_N=\xi_N\big)~=~\frac{1}{2^N}
\end{equation}
for any $\xi_1,\ldots,\xi_N\in\{-1,1\}$. The voting measure describes the mutual influence of the voters on each other, mathematically speaking it describes the correlation structure of the voting system. The above example describes \emph{independent} voters - in some sense the classical case of the theory. An extreme case is given by the measure $\IP_u$
\begin{align}\label{unanim}
\IP_u\big(X_1=1, X_2=1,\ldots, X_N=1\big)~&=~\IP_u\big(X_1=-1, X_2=-1,\ldots, X_N=-1\big)~\notag\\
&=~\frac{1}{2}\,.
\end{align}
For this (rather boring) voting measure the only possible outcomes are the unanimous votes, it represents total (positive) correlation.

If $\IP$ is a voting measure, we denote the expectation value with respect to $\IP$ by $\IE$, as was already anticipated in \eqref{expDelta}.
Since we assume that the numbers $N_\nu$ are odd, it follows that $S_\nu\not=0$. From this we conclude that $\IE(\chi_\nu)=0$ for any voting
measure.

\section{Optimal Weights for Independent States}\label{sec:optweight}
We begin by determining optimal weights, under the assumption that voters in different states are independent.
Thus, we assume that the random variables $\Xn$ and $X_{\kappa\;\!j}$ are independent for $\nu\not=\kappa$.

We want to minimize the function

\begin{align}
\ID(\gamma_1,\ldots,\gamma_M)~&=~\IE\big(\,\Delta(\gamma_1,\ldots, \gamma_M)^2\big)\notag\\
&=~\sum_{\nu,\kappa=1}^M\,\Big(\gamma_\nu \gamma_\kappa\IE\big(\chi_\nu \chi_\kappa\big)- 2\gamma_\nu\IE\big(\chi_\nu S_\kappa\big) + \IE\big(S_\nu S_\kappa\big)\Big)\,.
\end{align}
The function $\ID(\gamma_1,\ldots,\gamma_M)$ is a measure for the expected democracy deficit for voting weights $\gamma_1,\ldots,\gamma_M$.

By the assumption of independent \emph{states} we can conclude that
\begin{align}
\IE\big(\chi_\nu \chi_\kappa\big)~&=~\IE\big(\chi_\nu\big)\,\IE\big(\chi_\kappa\big)~=~0\qquad\tn{for $\nu\not=\kappa$}\,,\label{f1}\\
\IE\big(\chi_\nu S_\kappa\big)~&=~\IE\big(\chi_\nu\big)\,\IE\big(S_\kappa\big)~=~0\qquad\tn{for $\nu\not=\kappa$}\,,\\
\intertext{and}
\IE\big(S_\nu S_\kappa\big)~&=~\IE\big(S_\nu\big)\,\IE\big(S_\kappa\big)~=~0\qquad\tn{for $\nu\not=\kappa$}\,.
\end{align}
Moreover, we have $\chi_\nu^2=1$ and $\chi_\nu S_\nu=|S_\nu|$, thus
\begin{equation}
\ID(\gamma_1,\ldots,\gamma_M)~=~\sum_{\nu=1}^M\;\Big(\gamma_\nu^2\;-\; 2 \gamma_\nu \IE\big(|S_\nu|\big)\;+\;\IE\big(S_\nu^2\big)\Big)\,.\label{f2}
\end{equation}
It is not hard to find the minimizing weights $g_\nu$ (by the usual procedure: Find the zeros of the derivative),
in fact: The weights $g_1,\ldots,g_M$ which minimize the function $\ID$ are given by
\begin{equation}
g_\nu~=~\IE\big(\big|S_\nu\big|\big)\,.
\end{equation}
This result has a very intuitive interpretation. The quantity $S_\nu$ is the difference between the `yes'-votes and the `no'-votes,
so $|S_\nu|$ describes the margin of the voting outcome, i.~e., the surplus of votes of the winning party. Therefore, the optimal
weights $g_\nu$ for the state $\nu$ are given by the expected margin of a vote in that state. In fact, the delegate of state $\nu$
does not represent the opinion of \emph{all} voters in this state, but only those who agree with the majority, he or she acts against
the will of the minority, so as a net result the delegate just represents the margin.

We can also compute the expected democracy deficit $\ID$ for the optimal weights $g_1,\ldots,g_M$
\begin{equation}
\ID(g_1,\ldots,g_M)~=~\sum_{\nu=1}^M\;\Big(\IE\big(\big|S_\nu\big|^2\big)-\IE\big(\big|S_\nu\big|\big)^2\Big)~=~\sum_{\nu=1}^M\;\IV\big(\big|S_\nu\big|\big)
\end{equation}
where $\IV(|S_\nu|)$ denotes the variance of the random quantity $|S_\nu|$.

We emphasize that we did not yet make assumptions about the correlation structure of voters \emph{inside} a country. Of course, the numerical evaluation
of the optimal weights and minimal democracy deficit requires further assumptions on the correlation between voters.

\section{Independent Voters}\label{sec:indep}

In this section we assume that all voters act independently of each other, in mathematical terms: all random variables $\Xn$ are independent of each other.
Under this assumption we can compute the optimal weight $g_\nu=\IE(|S_\nu|)$ as well as the minimal expected democracy deficit.

For the independent random variables $\Xn$  we have the central limit theorem, namely the weighted sums
\begin{equation}
\frac{1}{\sqrt{N_\nu}}\,S_\nu~:=~\frac{1}{\sqrt{N_\nu}}\;\sum_{i=1}^{N_\nu}\;\Xn
\end{equation}
are asymptotically distributed for large $N_\nu$ according to a standard normal distribution (cf. Lamperti \cite{Lamperti1996}). From this it follows that for large $N_\nu$
\begin{align}
\IE\big(\big|S_\nu\big|\big)~&\approx~\frac{\sqrt{2}}{\sqrt{\pi}}\,\sqrt{N_\nu}\,,\\
\IE\big(\big|S_\nu\big|^2\big)~&\approx~\,\sqrt{N_\nu}\,,\\
\intertext{and}
\IV\big(\big|S_\nu\big|\big)~&\approx~\frac{\pi-2}{\pi}\,N_\nu\,.
\end{align}
We conclude that the optimal weight for independent voters is proportional to the square root of the population. This is exactly the content of
the square root law by Penrose (see \cite{Penrose1946} and \cite{FeMa1998}).

The above formulae also allow us to evaluate the minimum of the expected democracy deficit
\begin{equation}
\ID(g_1,\dots,g_M)~\approx~\frac{\pi-2}{\pi}\,N\,.
\end{equation}
This implies that the \emph{expected democracy deficit per voter}, namely
\begin{equation}
\IE\left(\left(\frac{\Delta}{N}\right)^2\right)
\end{equation}
converges to zero as $N$ becomes large (with convergence rate $\tfrac{1}{N})$.

\section{The Collective Bias Model}

Now, we introduce and discuss a model for collective behaviour of voters. The basic idea is that there is a
mainstream opinion, e.~g., a common belief due to the country's tradition or the influence of opinion makers.
For a given proposal $\omega$ we model this `common belief' by a value $\zeta\in[-1,1]$ which depends on the proposal at hand.
The value $\zeta=1$ means there is such a strong common belief in favor of the proposal that all voters will vote `yes',  $\zeta=-1$
means all voters will vote `no'. In general, $\zeta$ denotes the expected outcome of the voting, i.~e., $\IE(\Xn)$. The voting results
$\Xn$ themselves fluctuate around this value randomly.

Let us be more precise about this. Suppose the voting results are $X_1,\ldots,X_N$ (where we dropped the index $\nu$ for notational simplicity).
Let $\mu$ be a measure on $[-1,1]$, which is the distribution of the common belief value $\zeta$, that is
$\mu(]a,b[)$ is the probability that the value $\zeta$ is between $a$ and $b$. Let $P_\zeta$ be the probability measure on $\{-1,1\}$
with
\begin{align*}
P_\zeta(X_1=1)=p_\zeta=\tfrac{1}{2}(1+\zeta)\,,
\end{align*} so that
\begin{align*}
E_\zeta(X_1):=P_\zeta(X_1=1)-P_\zeta(X_1=-1)=p_\zeta-(1-p_\zeta)=\zeta\,.
\end{align*}
For a given value of $\zeta$ we set
\begin{equation}
\P_\zeta(\xi_1,\ldots,\xi_N)~=~\prod_{i=1}^N\;P_\zeta(\xi_i)\,.
\end{equation}
For any $\zeta\in[-1,1]$ the expression $\P_\zeta$ is a probability distribution on $\Omega=\{-1,1\}^N$.
We define the collective bias measure $\IP_\mu$ with respect to $\mu$ as
\begin{equation}
\IP_\mu(X_1=\xi_1,\ldots,X_N=\xi_N)~:=~\int\;\P_\zeta(\xi_1,\ldots,\xi_n)\;d\mu(\zeta)\,.
\end{equation}

Note, that $\P_\zeta$ is \emph{not} a voting measure (unless $\zeta=\tfrac{1}{2}$). However $\IP_\mu$ \emph{is} a voting measure
if $\mu$ is invariant under sign change, i.~e., $\mu(]a,b[)=\mu(]-b,-a[)$. We call $\mu$ the \emph{bias measure}.

If the measure $\mu$ is concentrated in $0$, then $\IP_\mu$ makes the voting results $X_i$ independent, thus we are in the case of section \ref{sec:indep}.
If $\mu$ is the uniform distribution on $[-1,1]$ (that is every point is equally likely), then the corresponding measure was already considered by Straffin \cite{Straffin1977}
where he established an intimate connection of this model to the Shapley-Shubik index. In a similar way, the Penrose-Banzhaf measure is connected with the model of
independent voters.

The CBM can be looked upon as a model for spins in statistical mechanics. There the voters are replaced with elementary magnets (spins) which can be directed upwards ($X_i=1$)
or downwards ($X_i=-1$). In this language the Collective Bias Model describes spins which do not interact with each other but are influenced by an exterior magnetic field,
namely the collective bias $\zeta$.

In the papers \cite{Kirsch2007}, \cite{KiLa2012} and \cite{Langner2012} we investigate also another model for collective voting behaviour which comes directly from
statistical physics, the Curie-Weiss Model (CWM). In this model the spins (voters) influence each other by an interaction which makes spins to prefer to be directed
parallel to the others. For voting this means that voters prefer to agree to the other voters. The Curie-Weiss Model is a very interesting tool to investigate collective
behaviour. However, it is technically more involved than the other models we discuss. Therefore, we will mention it only rather briefly and refer to the papers mentioned above
for more details.

Let us define
\begin{equation}
H(X_1,\ldots,X_N)~=~-\frac{1}{N}\Big(\sum_{i=1}^N\,X_i\Big)^2\,.
\end{equation}
This is the \emph{energy function} for the spin configuration $X_1,\ldots, X_N$. We use this to define measures
\begin{equation}
    Q_\beta(X_1,\ldots, X_N)~=~e^{-\beta H(X_1,\ldots, X_N)}
\end{equation}
where $\beta\in]0,\infty[$ is the inverse temperature in statistical physics.
As a rule, $Q_\beta$ is not a probability measure, so we normalize it by dividing through its total mass $Z$ and set
\begin{equation}\label{CW}
    P_\beta(X_1,\ldots, X_N)~=~\frac{e^{-\beta H(X_1,\ldots, X_N)}}{Z}\,.
\end{equation}
This is the Curie-Weiss measure for inverse temperature $\beta$. The parameter $\beta$ measures the strength of the
interaction between the voters. The extreme case $\beta=0$ corresponds to the model of independent voters, the other
extreme $\beta=\infty$ describes the case of the measure $\IP_u$ defined in \eqref{unanim} for unanimous voting.

\section{Optimal weights for the Collective Bias Model}
Let us now suppose that voters in different countries are independent, but voting inside the countries follows the CBM with bias measure $\mu$.
According to section \ref{sec:optweight} in this case the optimal weights  are given by
\begin{equation}
g_\nu~=~\IE_\mu\big(\big|S_\nu\big| \big)\,.
\end{equation}
For large $N_\nu$ we have
 \begin{equation}\label{optwCBM}
g_\nu~=~\IE_\mu\big(\big|S_\nu\big| \big)~=\mu_1\,N_\nu
\end{equation}
where $\mu_1=\int|\zeta|d\mu(\zeta)$ is the first absolute moment of $\mu$. Note, that for any probability measure $\mu$ the quantity $\mu_1$ is non zero,
except for the case $\mu=\delta_0$, the measure is concentrated at the point $0$. This means that the optimal weights for a council are proportional to
the population of the respective country if the voters can be described by a CBM. This also includes the Straffin case ($\mu$ is the uniform distribution),
which corresponds to the Shapley-Shubik power index.

The only exception from proportionality is the case $\mu=\delta_0$ corresponding to independent voting (the Penrose-Banzhaf case), where the square root
law applies.

We mention that there is a `phase transition' for the Curie-Weiss Model if we vary $\beta$ from $0$ to $\infty$, namely

\begin{equation}
g_\nu~=~\IE_\beta\big(\big|S_\nu\big|\big)~=~
\left\{
  \begin{array}{ll}
    \frac{\sqrt{2}}{\sqrt{\pi}\sqrt{1-\beta}}\,\sqrt{N_\nu}, & \hbox{for $\beta<1$;} \\
\\
    C\,N_\nu^{\tfrac{3}{4}}, & \hbox{for $\beta=1$;} \\
\\
    C(\beta)\,N_\nu, & \hbox{for $\beta>1$.}
  \end{array}
\right.
\end{equation}
The constant $C(\beta)$ converges to $0$ as $\beta\searrow 1$ and to $1$ as $\beta\nearrow\infty$.

\section{Democracy Deficit for the Collective Bias Model}
Given the optimal weights \eqref{optwCBM} for the CBM (and independent states) we can compute (the asymptotic behaviour of) the expected democracy Deficit $\ID_\mu$

\begin{equation}
\ID_\mu~=~\sum_{\nu=1}^M\;\IV\big(\big|S_\nu\big|\big)~\approx~(\mu_2-\mu_1^2)\,N^2
\end{equation}
where $\mu_1=\int|\zeta| d\mu(\zeta)$ and $\mu_2=\int|\zeta|^2 d\mu(\zeta)$. Note that $\mu_2-\mu_1^2\not=0$ unless $\mu$ is concentrated in at most two points.
It follows that the expected democracy deficit per voter, i.~e.,
\begin{align*}
\IE_\mu\left(\left(\frac{\Delta}{N}\right)^2\right)
\end{align*}
converges to a \emph{positive} constant as the $N_\nu$ tend to infinity (in a uniform way, i.~e., $N_\nu=\alpha_\nu N$).

It is interesting to remark that the expected democracy deficit per voter converges also to a constant if we choose a \emph{non optimal} voting weight, like
for instance $g_\nu\sim\sqrt{N_\nu}$ or $g_\nu=1$ for all $\nu$.
This constant will in general be larger than the one for the optimal weights, but the order of magnitude of $\ID$ is \emph{not} changed.

For the Curie-Weiss Model the expected democracy deficit per voter converges to zero (for $\beta\not=1$ even with rate $\tfrac{1}{N}$).

\section{A Model with Global Collective Behaviour}

So far we have always assumed that voter in different states act independently. In this section we consider the case of collective behaviour across country borders.
We assume that \emph{all} voters act according to the Collective Bias measure $\IP_\mu$. This means there is a common belief, expressed through the measure $\mu$,
for all voters in the union.

Then, the formulae \eqref{f1} -- \eqref{f2} are no longer valid. In fact, determining the optimal voting weights requires to solve a rather complicated system of $M$ dependent linear
equations. Instead of doing this we try to look at the democracy deficit directly. It turns out that for large $N_\nu$ we have for any $\nu,\kappa$
\begin{align}
\IE_\mu\big(\chi_\nu \chi_\kappa\big)~&\approx 1\,,\label{f3}\\
\IE_\mu\big(\chi_\nu S_\kappa\big)~&\approx~\IE_\mu\big(|S_\kappa\big|)~\approx~\mu_1\,N_\kappa\,\\
\intertext{and}
\IE\big(S_\nu S_\kappa\big)~&\approx~\mu_2\;N_\nu\,N_\kappa\,.
\end{align}

Inserting these terms into the expression for $\ID$ we obtain
\begin{align}
\ID(g_1,\ldots,g_M)~\notag=~&\sum_{\nu,\kappa=1}^M\,\IE_\mu(\chi_\nu \chi_\kappa)\,g_\nu g_\kappa\;-\;2\sum_{\nu=1}^M\,g_\nu\,\sum_{\kappa=1}^M\,\IE_\mu(\chi_\nu S_\kappa)
\;+\;\sum_{\nu,\kappa=1}^M\,\IE_\mu(S_\nu S_\kappa)\notag\\
\approx~&\sum_{\nu,\kappa=1}^M\,g_\nu g_\kappa\;-\;2\sum_{\nu=1}^M\,g_\nu\,\sum_{\kappa=1}^M\,\mu_1\,N_\kappa
\;+\;\sum_{\nu,\kappa=1}^M\,\mu_2\,N_\nu N_\kappa\notag\\
=~&\left(\sum_{\nu=1}^M g_\nu\right)^2\;-\;2\mu_1\,\big(\sum_{\nu=1}^M g_\nu\big)\,N\;+\;\mu_2\,N^2\notag\\
=~&G^2\;-\;2\mu_1\,G\;+\;\mu_2\,N^2\,.\label{f4}
\end{align}
This last expression depends only on the sum $G=\sum_{\nu=1}^M g_\nu$ of the voting weights and \emph{not} on the single weight $g_\nu$.
This means that for large $N_\nu$ the asymptotic value of $\ID$ does not depend on the way the weights are distributed among the member states of the union.
The minimal value of $\ID$ is obtained by choosing $G=\mu_1 N$ independent of the values of the particular weight $g_\nu$. We also note that the value
of $G$ has no real meaning, since we don't change the voting system at all if we multiply all weights (and the quota) with the same number $C>0$.

Finally, we remark that the somewhat hand waving arguments in \eqref{f4} need a careful mathematical interpretation. A precise formulation gives:

\begin{equation}
\lim_{N\to\infty} \IE_\mu\Big(\Big(\frac{\Delta(g_1,\ldots,g_M)}{N}\Big)^2\Big)~=~\mu_2-\mu_1^2
\end{equation}
for $G=\sum_{\nu=1}^M g_\nu=\mu_1 N$ and
\begin{equation}
\liminf_{N\to\infty} \IE_\mu\Big(\Big(\frac{\Delta(g_1,\ldots,g_M)}{N}\Big)^2\Big)~\geq~\mu_2-\mu_1^2
\end{equation}
for any arbitrary choice of $g_\nu$.
This result can be interpreted in the following way: If there is a strong common belief in the union across border lines then it doesn't matter
how one distributes the voting weights in the council.

\end{document}